# Congruence Properties of Mersenne Primes

M. S. Srinath[1], Garimella Rama Murthy[2], V. Chandrasekaran[3]

[1,3] Sri Sathya Sai Institute of Higher Learning,
Puttaparthi, Anantapur, AP, INDIA

[2] Associate Professor,
International Institute of Information Technology (IIIT),
Gachibowli, Hyderabad-32, AP, INDIA

**ABSTRACT**

In this research paper, relationship between every Mersenne prime and certain Natural numbers is explored. We begin by proving that every Mersenne prime is of the form $\{4n + 3, \text{ for some integer 'n'}\}$ and generalize the result to all powers of 2. We also tabulate and show their relationship with other whole numbers up to 10. A number of minor results are also proved. Based on these results, approaches to determine the cardinality of Mersenne primes are discussed.

**1. Introduction:**

With the advent of the concept of "counting", homo-sapien civilization in various parts of the planet began a systematic study of the properties of whole numbers ( integers ). Prime numbers were identified as those whole numbers which have the numbers "one" and "itself" as the divisors. Mathematicians became interested in determining the cardinality of prime numbers. Euclid provided an elegant proof that the prime numbers are infinite ("countably") (Har). Also, the fundamental theorem of arithmetic was proved by the Greek mathematicians. It states that every integer can be written as the product of powers of finitely many primes (below the integer). Euclid also provided an algorithm for computing the Greatest Common Divisor (GCD) of two whole numbers/integers. Also few books in Euclid's Elements were devoted to number-theoretic investigations. Diophantus formulated the problem of solving algebraic equations in integers. For instance, it was proved that the sides of a right angle triangle satisfy an interesting quadratic Diophantine equation ( named as the Pythagorean Theorem ).

For quite a few centuries, there were only scattered mathematical results in number theory. Modern number theory began with the efforts of Fermat who recorded several interesting theorems ( related to whole numbers ) in the margins of his copy of "Arithmetica" ( written by Diophantus ) (Dun). For instance, the following theorem is the famous little theorem of Fermat:

**Theorem 1:** If "p" is a prime and "a" is a whole number which does not have "p" as a factor, then "p" divides evenly into $a^{p-1} - 1$.
Equivalently, in the number theoretic notation we have
$$a^{p-1} \equiv 1 \ (mod \ p)$$

Euler provided proofs for several theorems recorded by Fermat in his copy of "Arithmetica". As time progressed, generations of mathematicians became interested in the following (so called) arithmetical function (kan):

$$\pi(n) : \text{Total number of prime numbers until } 'n'$$

Gauss conjectured the so called "Prime Number Theorem" on empirical evidence. The formal proof was given by Hadamard et.al. using the mathematical tools of analytic number theory. The Theorem states that

$$\lim_{n \to \infty} \frac{\pi(n)}{\frac{n}{\log n}} = 1.$$

Ever since, Euclid provided the elegant proof of the infinitude of the primes, mathematicians were interested in identifying special structured classes of primes and their cardinality. For instance, Fermat conjectured that

$$F_n = 2^{2^n} + 1$$

is a prime number for all $n \geq 1$. Euler refuted the conjecture by proving that $F_5$ is composite. In fact, so far, for no other value of 'n' ( $\geq 5$ ), $F_n$ was found to be a prime number.

Mathematicians became interested in the problem of determining the cardinality of primes in arithmetical progressions like { 4 n + 1, $n \geq 1$ } and { 4 m + 3, $m \geq 1$ }. It was shown that Euclid's proof argument can be naturally extended to prove that there are infinitely many primes in the arithmetical progression { 4 n + 1, $n \geq 1$ }.

This approach did not directly work for determining the cardinality of primes in the arithmetical progression of the form { 4 m + 3, $m \geq 1$ }. Using the tools of analytic number theory, Dirichlet showed that there are infinitely many primes in every arithmetic progression of the form
{ a n + b, $n \geq 1$ } with (a,b) =1 i.e. "a" and "b" are relatively prime.

French father Marin Mersenne, defined primes of the form

$$M_p = 2^p - 1 \text{ for } 'p' \text{ a prime number}.$$

An interesting open problem is to determine the cardinality of Mersenne primes. There are several such open research problems. For instance, determination of cardinality of primes in the sequence $\{ n^2 + 1, n \geq 1 \}$ is also an open research problem. An interesting general result is the following:

" No polynomial can generate primes for all integral values of the independent variable".

In this research paper, we attempt the problem of determining the cardinality of Mersenne primes using congruence properties satisfied by Mersenne primes. This technical report is organized as follows. In Section 2, we summarize the Dirichlet Theorem on arithmetical primes. In Section 3, we explore the relationship between Mersenne primes and arithmetical primes. We prove a few interesting results. In Section 4, our efforts at determining the cardinality of Mersennes primes ( and thus that of even perfect numbers ) are summarized. In Section 5, we conclude the research paper.

**2. Arithmetical Primes : Dirichlet's Theorem**:

In view of the infinitude of primes of the form { 4 n + 1, $n \geq 1$ } , mathematicians were interested in the most general result regarding the cardinality of primes generated through an arithmetical progression. In this connection, the following observation follows readily:

- Suppose the Greatest Common Divisor (G.C.D) of a, b be D >1. Then it is easy to see that the arithmetic progression
{a n + b, $n \geq 1$} will contain no prime numbers. Thus, let us consider the non-trivial case i.e
G.C.D. of a,b = (a,b) =1 i.e. a,b are relatively prime.

Using the tools of analytic number theory, Dirichlet proved the following interesting Theorem:

**Theorem 2:** Suppose 'a', 'b' are relatively prime numbers. Then the arithmetical progression { a n + b, $n \geq 1$ } will contain infinitely many prime numbers.

In view of the above Dirichlet's theorem, mathematicians are interested in determining the cardinality of prime numbers generated by a quadratic, cubic and higher degree polynomial ( with the independent variable assuming integer values ). The authors are currently investigating this problem.

3. **Relationship between Mersenne Primes and Arithmetic Primes:**

On numerical (empirical) evidence, the authors expected a relationship between Mersenne primes and arithmetical primes. Mersenne primes are related to the primes of the form { 4 n + 3, $n \geq 1$ }. The following Theorem is an interesting result in this direction:

**Theorem 3:** Every Mersenne prime, $M_p$ is a prime number of the form
4 n + 3 for some integer "n".

**Proof:** First we consider the integers
$$T_n = 2^n - 1 \text{ for all integer values of n.}$$

We show that $T_n = 4k + 3$ for some integral value of k. The proof approach is based on mathematical induction. It is easy to see that for n=2,
$$T_2 = 3 = 4(0) + 3. \text{ Thus } k = 0.$$

Also, for n=3, we have that
$$T_3 = 7 = 4(1) + 3. \text{ Thus } k = 1.$$
Further, for n=4, we have that
$$T_4 = 15 = 4(3) + 3. \text{ Thus } k = 3.$$

- Now, let us assume that
$$T_n = 4(k_0) + 3 \text{ for some integer } k_0.$$
It is easy to see that
$$T_{n+1} = 2^{n+1} - 1 = (2^n)2 - 1$$
$$= (4k_0 + 4)2 - 1 = 8k_0 + 7$$

$$= 4(2k_0 + 1) + 3 = 4k' + 3 \text{ for } k' = 2k_0 + 1.$$

Thus, by mathematical induction the claim holds true for n+1. Since, we are only interested in prime exponents, the claim still holds true. Thus, every Mersenne prime $M_p$ is a prime of the form 4n+3 for some integer n.   Q.E.D.

**Corollary:** No Mersenne prime can be written as the sum of squares of two Integers

**Proof:** The corollary follows from the celebrated Theorem of Fermat on the primes of the form $\{4n+3, n \geq 1\}$.                    Q.E.D.

**Note**: Not every prime of the form $\{4n+3$, for some integer 'n'$\}$ is a Mersenne prime. For example 19, 23.

- The above Theorem enables us to invoke results related to arithmetical primes of the form $\{4n+3, n \geq 1\}$ in association with the Mersenne primes. Details are discussed in the following section.

- As a generalization of Theorem 3, we expect the following result. It is trivial to see that $a^p - 1$ is never a prime number for $a \geq 3$. It will be an odd/even composite number when 'a' is even/odd number.

**Claim (i):** If 'p' is a prime number and 'a' is an even number, then $a^p \cong 0 \ (mod \ 4)$ and $a^p - 1 \equiv 3 \ (mod \ 4)$.

**Proof:**
Since 'a' is even, there exists a natural number 'k' such that, $a = 2k$.
Now, $a^p = 2^p k^p$.
Since $p \geq 2$, 4 divides evenly into $2^p$ and therefore 4 divides evenly into $a^p$
Hence we have, $a^p \cong 0 \ (mod \ 4)$ and by subtracting 1 from both sides, we get
$$a^p - 1 \cong 3 \ (mod \ 4)$$                    Q.E.D.

**Claim (ii):** If 'p' is an odd prime number and 'a' is an odd number, then $a^{p-1} \equiv 1 \ (mod \ 4)$ and $a^p - 1 \equiv a - 1 \ (mod \ 4)$.

**Proof:**
If 'a' is odd, it can be easily seen that $a^2 \cong 1 \ (mod \ 4)$.
Therefore, for any even number 'n', we have $a^n \cong 1 \ (mod \ 4)$.
Furthermore, for every odd number 'm', we get $a^m \cong a (mod \ 4)$.
If 'p' is an odd prime, then $a^p - 1 \cong a - 1 \ (mod \ 4)$.
Also, $p - 1$ is even which implies, $a^{p-1} \equiv 1 \ (mod \ 4)$.        Q.E.D.

We have seen relationship of Mersenne numbers with 4. Now we can move a step further and ask the question "what is the relationship of Mersenne numbers with other whole numbers?" To understand the problem more clearly let us look into the following table:

| Primes ($p$) | 2 | 3 | 5 | 7 | 11 | 13 | 17 | 19 |
|---|---|---|---|---|---|---|---|---|
| $p(mod\ 3)$ | 2 | 0 | 2 | 1 | 2 | 1 | 2 | 1 |
| $p(mod\ 4)$ | 2 | 3 | 1 | 3 | 3 | 1 | 1 | 3 |
| | | | | | | | | |
| $M_p = 2^p - 1$ | 3 | 7 | 31 | 127 | 2047 | 8191 | 131071 | 524287 |
| Is $M_p$ Prime? (y/n) | y | y | y | y | n | y | y | y |
| $M_p(mod\ 2)$ | 1 | 1 | 1 | 1 | 1 | 1 | 1 | 1 |
| $M_p(mod\ 3)$ | 0 | 1 | 1 | 1 | 1 | 1 | 1 | 1 |
| $M_p(mod\ 4)$ | 3 | 3 | 3 | 3 | 3 | 3 | 3 | 3 |
| $M_p(mod\ 5)$ | 3 | 2 | 1 | 2 | 2 | 1 | 1 | 2 |
| $M_p(mod\ 6)$ | 3 | 1 | 1 | 1 | 1 | 1 | 1 | 1 |
| $M_p(mod\ 7)$ | 3 | 0 | 3 | 1 | 3 | 1 | 3 | 1 |
| $M_p(mod\ 8)$ | 3 | 7 | 7 | 7 | 7 | 7 | 7 | 7 |
| $M_p(mod\ 9)$ | 3 | 7 | 4 | 1 | 4 | 1 | 4 | 1 |
| $M_p(mod\ 10)$ | 3 | 7 | 1 | 7 | 7 | 1 | 1 | 7 |

The above table contains primes up to 19, the corresponding Mersenne numbers and their congruence with whole numbers up to 10. This table encourages us to state the following results:

**Theorem 4:** Let $p > 2$ be an odd prime and $M_p$ be the corresponding Mersenne number, then, the following statements hold:
1. $M_p \cong 1(mod\ 2)$
2. $M_p \cong 1(mod\ 3)$
3. $M_p \cong 3(mod\ 4)$
4. $M_p \cong 1(mod\ 5)$ if $p \cong 1(mod\ 4)$
   $M_p \cong 2(mod\ 5)$ if $p \cong 3(mod\ 4)$
5. $M_p \cong 1(mod\ 6)$
6. $M_p \cong 1(mod\ 7)$ if $p \cong 1(mod\ 3)$, $p > 3$
   $M_p \cong 3(mod\ 7)$ if $p \cong 2(mod\ 3)$
7. $M_p \cong 7(mod\ 8)$
8. $M_p \cong 1(mod\ 9)$ if $p \cong 1(mod\ 3)$, $p > 3$
   $M_p \cong 4(mod\ 9)$ if $p \cong 2(mod\ 3)$
9. $M_p \cong 1(mod\ 10)$ if $p \cong 1(mod\ 4)$
   $M_p \cong 7(mod\ 10)$ if $p \cong 3(mod\ 4)$

The above results can be proved in the similar lines as proof of Theorem 3. Note that the above results suggest that Mersenne numbers are not divisible by any of the numbers below 10.

Interestingly we don't have to find the congruence of a Mersenne number with every whole number. It is sufficient to know their congruence with primes and powers of primes. Note by Fundamental Theorem of Arithmetic, any whole number can be written as product of powers of primes. So that by invoking Chinese Remainder Theorem (CRT), the congruence can be known.

For example: To know the value of $M_p(mod\ 12)$ we combine $M_p \cong 1(mod\ 3)$ and $M_p \cong 3(mod\ 4)$ using CRT. After solving, we get the solution to be 7. Hence,
$$M_p \cong 7(mod\ 12)$$

We now state a generalization of Theorem 3 to the powers of 2. Let us consider the powers of 2; $2^i$ where $i$ is a Natural number.

We would like to know $M_p \pmod{2^i}$.

If $i \geq p$, then $2^i > 2^p - 1$, therefore $M_p \pmod{2^i} = M_p$.

Suppose, $i < p$, then, $2^i$ divides $2^p$. Therefore, $2^p \cong 0 \pmod{2^i}$
$$\Rightarrow M_p \cong -1 \pmod{2^i}$$

Hence the following theorem:

**Theorem 5:** Let p be a prime and $M_p$ be the corresponding Mersenne number, then
$$M_p \cong -1 \pmod{2^i} \quad \forall \, i < p$$

We now state and prove a theorem on the relation of $M_p$ with p itself.

**Theorem 6:** Let p be an odd prime and $M_p$ be the corresponding Mersenne number, then $M_p \cong 1 \pmod{p}$

**Proof:** Since p is an odd prime $gcd(2, p) = 1$. Hence by invoking Fermat's Little Theorem, we have,
$$2^{p-1} \cong 1 \pmod{p}$$
$$\Rightarrow 2^p \cong 2 \pmod{p}$$
$$2^p - 1 \cong 1 \pmod{p}$$
$$\Rightarrow M_p \cong 1 \pmod{p}$$

We can also arrive at some more results by combining Theorems 4, 5 and 6 using Chinese Remainder Theorem (CRT).

For example, combining Theorem 4 (part 5) with Theorem 6 using CRT yields,
$$M_p \cong 1 \pmod{6p} \quad \forall \, p > 3$$

4. **Cardinality of Mersenne Primes and Perfect Numbers**:

For a long time, determination of cardinality ( finiteness or countably infiniteness ) of Mersenne primes has been an open problem. It has been well known that this problem is directly related to the determination of cardinality of even perfect numbers ( in view of the result of Euler on the structure of even perfect numbers ). Specifically, we have the following Theorem:

**Theorem 7:** An even number, N is a perfect number if and only if

$$N = 2^{p-1} M_p, \text{ where}$$

'p' is a prime number and $M_p$ is the associated Mersenne prime number.

- Another interesting problem is to determine whether there are any odd perfect numbers or not. Equivalently, it is to be determined whether all perfect numbers are just even numbers only.

Based on the Theorem 3, we attempt the problem of determination of cardinality of Mersenne primes. The following definition is utilized.

**Definition:** A Diophantine equation is called a "primal Diophantine equation" if all the integral solutions are required to be the prime numbers.

Based on the Theorem 3, we need to find the cardinality of solutions of the following "primal" Diophantine equation:
$$M_p = 2^p - 1 = q = 4k + 3,$$
where 'p', 'q' are prime numbers and 'k' is an integer.
i.e The problem effectively boils down to the answer to the following question:

Q: Are there infinitely many arithmetic primes of the form $\{4k + 3, k \geq 1\}$ that are also Mersenne primes?

- The following Lemma directly follows as a special case of Fermat's Little Theorem:

**Lemma 1:** If 'p' is a prime other than 2 and $M_p = 2^p - 1$ is the corresponding Mersenne prime, then we have
$$M_p = 2pk + 1 \quad \text{for some integer } 'k'.$$

- The following Lemmas / Theorems readily follow from the approach to prove the Fermat's Little Theorem.

**Theorem 8:** If 'p' is a prime and 'a' is any whole number, then 'p' divides evenly into $a^p - a$.

**Theorem 9:** If 'p' is a prime and 'a' is a whole number, then $(a+1)^p - (a+1)$ is evenly divisible by 'p'.

**Theorem 10:** If 'p' is a prime and 'a' is a whole number, then
$(a+1)^p - (a^p + 1)$ is evenly divisible by 'p'.
- Using Theorem 5 with $a = 2$, we readily have that
$$M_p = pk + 1 \quad \text{for some integer } 'k'.$$

- Special cases of Theorems 9 and 10 with $a = 2$ can easily be derived. They are avoided here for brevity.

Theorem 3 and Lemma 1 can be combined to form a stronger necessary condition for Mersenne primes. Theorem 3 states, "Every Mersenne prime is of the form $4n + 3$ for some integer 'n'." And Lemma 1 states, "Mersenne prime $M_p$ ($p > 2$) is of the form $2pk + 1$ for some integer 'k'. Combining these two we get,
If $M_p$ ($p > 2$) is a Mersenne prime, then there exists integers 'k' and 'n' such that
$$M_p = 2pk + 1 = 4n + 3$$
By solving, we get $pk = 2n + 1$     --- (1)
RHS of the above equation (1) is an odd number, 'p' is an odd prime which forces 'k' to be odd.
Hence we can state the following:

**Theorem 11:** If 'p' is an odd prime and $M_p = 2^p - 1$ is the corresponding Mersenne prime, then there exists an odd number 'k' such that $M_p = 2pk + 1$.

Since 'k' is odd, we can replace 'k' by $2k + 1$ and rephrase the above statement as:
**Theorem 12:** If 'p' is an odd prime and $M_p = 2^p - 1$ is the corresponding Mersenne prime, then there exists an integer 'k' such that $M_p = 4pk + 2p + 1$.

Equation (1) also means that 'p' divides evenly into $2n + 1$. So, We can also state the following:

**Theorem 13:** If 'p' is an odd prime and $M_p = 2^p - 1$ is the corresponding Mersenne prime, then there exists an integer 'n' with the property $p|(2n + 1)$ such that
$$M_p = 4n + 3$$

So far in this section, we have only considered the sequence $\{4n + 3\}$. It should be noted that a similar analysis can be done using the sequences obtained by the virtue of Theorem 4, 5, 6 or a CRT combination of them.

Some of the worthwhile sequences are $\{3n + 1\}, \{6n + 1\}, \{8n + 7\}, \{12n + 7\}$, $\{pn + 1\}$ and $\{6pn + 1\}$

The following approaches are currently attempted to determine the cardinality of Mersenne primes:

**Approach 1:** Utilize the essential idea of Euclid's proof and the fact that the product of finitely many Mersenne primes is a product of finitely many geometric progressions

**Approach 2:** Again utilize the idea of Euclid's proof along with the following fact:
$$M_p = (1+1)^p - 1 = \binom{p}{1} + \binom{p}{2} + \cdots + \binom{p}{p} = 1 + 2p + \binom{p}{2} + \cdots + \binom{p}{p-2}.$$

**Approach 3:** Utilize the powerful tools of analytic number theory (as in the case of proof of Dirichlet's Theorem) alongside Theorem 3, 4, 5 and 6 proved in this research paper.

**5. Conclusions:**

In this research paper, it is proved that every Mersenne prime is a prime of the form { 4 n + 3, for some integer 'n' }. Necessary conditions on Mersenne primes in terms of congruences modulo certain natural numbers are also proved. Generalizations of these results using CRT are briefly discussed. Some special cases of the Fermat's Little Theorem are briefly considered. Approaches to establish the cardinality of Mersenne primes are summarized.